\documentclass[12pt,reqno]{amsart}
\usepackage{enumerate, latexsym, amsmath, amsfonts, amssymb, amsthm, color}
\def\pmod #1{\ ({\rm{mod}}\ #1)}
\def\Z{\Bbb Z}

\def\Q{\Bbb Q}

\def\C{\Bbb C}

\def\l{\left}
\def\r{\right}
\def\bg{\bigg}
\def\({\bg(}
\def\){\bg)}
\def\t{\text}
\def\f{\frac}
\def\mo{{\rm{mod}\ }}

\def\ls{\leqslant}

\def\sm{\setminus}

\def\ve{\varepsilon}

\def\eq{\equiv}

\def\da{\delta}

\def\Proof{\noindent{\it Proof}}
\def\Ack{\medskip\noindent {\bf Acknowledgment}}
\theoremstyle{plain}
\newtheorem{theorem}{Theorem}

\newtheorem{lemma}{Lemma}
\newtheorem{corollary}{Corollary}
\newtheorem{conjecture}{Conjecture}
\theoremstyle{definition}

\theoremstyle{remark}
\newtheorem{remark}{Remark}

\makeatletter
\@namedef{subjclassname@2020}{%
  \textup{2020} Mathematics Subject Classification}
\makeatother
 \vspace{4mm}

\begin{document}

\hbox{Accepted by Publ. Math. Debrecen}
\medskip

\title
[{Trigonometric identities and quadratic residues}]
{Trigonometric identities\\ and quadratic residues}

\author
[Zhi-Wei Sun] {Zhi-Wei Sun}

\address{Department of Mathematics, Nanjing
University, Nanjing 210093, People's Republic of China}
\email{zwsun@nju.edu.cn}

\keywords{Trigonometric functions, combinatorial identities, quadratic residues, roots of unity.
\newline \indent Supported by the Natural Science Foundation of China (grant no. 11971222).}

\begin{abstract}  In this paper we obtain some novel identities involving trigonometric functions.
Let $n$ be any positive odd integer. We mainly show that
$$\sum_{r=0}^{n-1}\frac1{1+\sin2\pi\frac{x+r}n+\cos2\pi\frac{x+r}n}
=\frac{(-1)^{(n-1)/2}n}{1+(-1)^{(n-1)/2}\sin 2\pi x+\cos 2\pi x}$$
for any complex number with $x+1/2,x+(-1)^{(n-1)/2}/4\not\in\Z$, and
$$\sum_{j,k=0}^{n-1}\frac1{\sin 2\pi\frac{x+j}n+\sin2\pi \frac{y+k}n}=\frac{(-1)^{(n-1)/2}n^2}{\sin 2\pi x+\sin2\pi y}$$
for all complex numbers $x$ and $y$ with $x+y,x-y-1/2\not\in\Z$. We also determine
the values of $\prod_{k=1}^{(p-1)/2}(1+\tan\pi\frac{k^2}p)$ and $\prod_{k=1}^{(p-1)/2}(1+\cot\pi\frac{k^2}p)$ for any odd prime $p$.
In addition, we pose several conjectures on the values of $\prod_{k=1}^{(p-1)/2}(x-e^{2\pi ik^2/p})$
with $p$ an odd prime and $x$ a root of unity.
\end{abstract}
\maketitle

\section{Introduction}
\setcounter{lemma}{0}
\setcounter{theorem}{0}
\setcounter{corollary}{0}
\setcounter{remark}{0}
\setcounter{equation}{0}

For the classical $\Gamma$-function, the Gauss multiplication formula states that
for each $n\in\Z^+=\{1,2,3,\ldots\}$
we have
$$\prod_{r=0}^{n-1}\Gamma\l(z+\f rn\r)=(2\pi)^{(n-1)/2}n^{1/2-nz}\Gamma(nz)$$
for all $z\in\C$ with $nz\not\in\{0,-1,-2,\ldots\}$, where $\C$ denotes the field of complex numbers.
As $\Gamma(z)\Gamma(1-z)=\pi/(\sin \pi z)$ for all $z\in\C\sm\Z$, Gauss' multiplication formula
implies the known formula
\begin{equation*}\prod_{r=0}^{n-1}\l(2\sin\pi\f{x+r}n\r)=2\sin \pi x\ \ \ (n\in\Z^+\ \t{and}\ x\in\C).
\end{equation*}
Taking logarithmic derivative for the last equality, one gets the known identity
$$\f1n\sum_{r=0}^{n-1}\cot\pi\f{x+r}n=\cot\pi x\ \ (n\in\Z^+\ \t{and}\ x\in\C\sm\Z).$$
Taking the derivatives of both sides of this formula, we obtain another well-known formula
\begin{equation}\label{csc2}\f1{n^2}\sum_{r=0}^{n-1}\csc^2\pi\f{x+r}n=\csc^2\pi x\ \ (n\in\Z^+\ \t{and}\ x\in\C\sm\Z).\end{equation}
If $n$ is a positive odd integer, then by taking $x=n/2$ in \eqref{csc2} we get the known identity
\begin{equation}\label{secant}\sum_{r=0}^{n-1}\sec^2\pi\f rn=n^2.
\end{equation}

There are lots of work on trigonometric power sums (cf. \cite{WZ94,BY,WZ}).
In 1994 K. S. Williams and N.-Y. Zhang \cite[Section 3]{WZ94}
expressed the sums
$$\sum_{r=1}^{n-1}\sin2\pi\f{ar}n\cot^{2m}\pi\f rn\ \ \t{and}\ \ \sum_{r=1}^{n-1}\cos2\pi\f{ar}n\cot^{2m}\pi\f rn
$$ with $a\in\{1,\ldots,n-1\}$ in terms of Bernoulli polynomials.
In 2002, B. C. Berndt and B. P. Yeap \cite[Corollary 2.2]{BY} gave a formula for $\sum_{r=1}^{n-1}\cot^{2m}\pi\f rn$ in terms of Bernoulli numbers.
For example, it is known (cf. \cite[Corollaries 2.3 and 2.6]{BY}) that
$$\sum_{r=1}^{n-1}\cot^2\pi\f rn=\f{(n-1)(n-2)}3$$
and $$\sum_{r=1}^{n-1}\cot^4\pi\f rn=\f{(n-1)(n-2)}{45}(n^2+3n-13).$$

In this paper we obtain some new trigonometric identities. Recall that for any positive odd integer $n$ we have
$$\l(\f{-1}n\r)=(-1)^{(n-1)/2}\ \ \t{and}\ \ \l(\f 2n\r)=(-1)^{(n^2-1)/8},$$
where $(\f{\cdot}n)$ denotes the Jacobi symbol.

Now we state our main results.

\begin{theorem}\label{Th1.1} Let $n$ be any positive odd integer. Then
\begin{equation}\label{sincos}\sum_{r=0}^{n-1}\f1{1+\sin2\pi\f{x+r}n+\cos2\pi\f{x+r}n}
=\f{(\f{-1}n)n}{1+(\f{-1}n)\sin 2\pi x+\cos 2\pi x}
\end{equation}
for any $x\in\C$ with $x+1/2,x+(-1)^{(n-1)/2}/4\not\in\Z$, and
\begin{equation}\label{-sincos}\sum_{r=0}^{n-1}\f1{1+\sin2\pi\f{x+r}n-\cos2\pi\f{x+r}n}
=\f{(\f{-1}n)n}{1+(\f{-1}n)\sin 2\pi x-\cos 2\pi x}
\end{equation}
for all $x\in\C$ with $x,x+(-1)^{(n-1)/2}/4\not\in\Z$.
\end{theorem}

\begin{corollary} We have
\begin{equation}\label{csc}\f1n\sum_{r=0}^{n-1}\csc 2\pi\f{x+r}n=\csc 2\pi x
\end{equation}
for all $x\in\C$ with $2x\not\in\Z$, and
\begin{equation}\label{sec}\f1n\sum_{r=0}^{n-1}\sec 2\pi\f{x+r}n=\l(\f{-1}n\r)\sec 2\pi x
\end{equation}
for any $x\in\C$ with $4x$ not an odd integer.
\end{corollary}
\begin{remark}\label{Rem1.1} For positive odd integers $m$ and $n$, X. Wang and D.-Y. Zheng \cite[p.\,1024]{WZ} expressed the sum
$\sum_{k=0}^{n-1}(-1)^k\sec^m \pi\f{x+k}n$ in terms of powers of $\sec \pi x$. We note that
\begin{align*}\sum_{k=0}^{n-1}(-1)^k\sec^m\pi\f{2x+k}n
=&\sum_{j=0}^{(n-1)/2}\sec^m2\pi\f{x+j}n
-\sum_{j=1}^{(n-1)/2}\sec^m\pi\f{2x+(n-2j)}n
\\=&\sum_{j=0}^{(n-1)/2}\sec^m2\pi\f{x+j}n+\sum_{j=1}^{(n-1)/2}\sec^m2\pi\f{x-j}n
\\=&\sum_{r=0}^{n-1}\sec^m2\pi\f{x+r}n.
\end{align*}
\end{remark}

\begin{corollary}\label{Cor1.2} For any positive odd integer $n$, we have
\begin{align}\label{sincos0}\sum_{r=0}^{n-1}\f1{1+\sin2\pi r/n+\cos2\pi r/n}=&\l(\f{-1}n\r)\f n2,
\\\label{-sincos0}\sum_{r=0}^{n-1}\f1{1+\sin\pi(2r+1)/n-\cos\pi(2r+1)/n}=&\l(\f{-1}n\r)\f n2,
\end{align}
\begin{equation}\label{sec0}\sum_{r=0}^{n-1}\sec2\pi\f rn=\l(\f{-1}n\r)n\ \t{and}\
\sum_{r=0}^{n-1}\sec\pi\f{2r+1}n=-\l(\f{-1}n\r)n.
\end{equation}
\end{corollary}

\begin{remark}\label{Rem1.2} \eqref{sincos} with $x=0$ and \eqref{-sincos} with $x=1/2$ yield \eqref{sincos0} and \eqref{-sincos0}. Putting $x=0,1/2$ in \eqref{sec} we get \eqref{sec0}.
Note also the simple trick
\begin{align*}&\sum_{r=0}^{n-1}\sec2\pi \f rn+\sum_{r=0}^{n-1}\sec\pi\f{2r+1}n
\\=&\sum_{k=0}^{2n-1}\sec\pi\f kn
=\sum_{k=0}^{n-1}\l(\sec\pi\f k{n}+\sec\pi\f{n+k}n\r)=0.
\end{align*}
\end{remark}

\begin{theorem}\label{Th1.2} Let $n$ be any positive odd integer.
Then
\begin{equation}\label{sin}\sum_{j,k=0}^{n-1}\f1{\sin 2\pi(x+j)/n+\sin2\pi (y+k)/n}=\l(\f{-1}n\r)\f{n^2}{\sin 2\pi x+\sin2\pi y}
\end{equation}
for all $x,y\in\C$ with $x+y\not\in\Z$ and $x-y-1/2\not\in\Z$, and
\begin{equation}\label{cos}
\sum_{j,k=0}^{n-1}\f1{\cos 2\pi(x+j)/n+\cos2\pi (y+k)/n}=\f{n^2}{\cos 2\pi x+\cos2\pi y}
\end{equation}
for all $x,y\in\C$ with $x\pm y-1/2\not\in\Z$. Also,
\begin{equation}\label{mix}\sum_{j,k=0}^{n-1}\f1{\sin 2\pi(x+j)/n+\cos2\pi (y+k)/n}=\f{n^2}{(\f{-1}n)\sin 2\pi x+\cos2\pi y}
\end{equation}
for all $x,y\in\C$ with $x\pm y+(-1)^{(n-1)/2}/4\not\in\Z$.
\end{theorem}

\begin{corollary}\label{Cor1.3} For any positive odd integer $n$, we have
\begin{align}\label{mix0}\sum_{j,k=0}^{n-1}\f1{\sin 2\pi j/n+\cos2\pi k/n}=&n^2
\\\label{mix1}\sum_{j,k=0}^{n-1}\f1{\sin\pi(2j+1)/n+\cos\pi(2k+1)/n}=&-n^2,
\\\label{cos0}\sum_{j,k=0}^{n-1}\f1{\cos 2\pi j/n+\cos2\pi k/n}=&\f{n^2}2,
\\\label{cos1}\sum_{j,k=0}^{n-1}\f1{\cos\pi(2j+1)/n+\cos\pi(2k+1)/n}=&-\f{n^2}2.
\end{align}
\end{corollary}
\begin{remark}\label{Rem1.3} The identities \eqref{mix0} and \eqref{mix1}
follow from \eqref{mix} with $x=y\in\{0,1/2\}$.
The identities \eqref{cos0} and \eqref{cos1} are just \eqref{cos} in the special case $x=y\in\{0,1/2\}$.
On August 2, 2019, the author posed \eqref{cos0} to {\tt MathOverflow} (cf. \cite{S19c}), and then both the user Wojowu and Fedor Petrov
provided proofs of \eqref{cos0}.
\end{remark}

\begin{corollary}\label{Cor1.4} Let $p$ be a prime with $p\eq3\pmod4$. Then
\begin{equation}\label{cosp}\sum_{1\ls j<k\ls(p-1)/2}\f1{\cos2\pi j^2/p+\cos2\pi k^2/p}=-\f{p+1}4\cdot\f{p-3}4.
\end{equation}
\end{corollary}

\begin{remark}\label{Rem1.4} Actually, the author found the identity in Corollary 1.4 inspired by his recent paper \cite{S19b}
on quadratic residues modulo primes, and this is the main motivation of this paper.
\end{remark}

Besides Theorems 1.1-1.2, we also
establish three other theorems on trigonometric identities related to
quadratic residues modulo primes.

\begin{theorem}\label{Th1.3} Let $p>3$ be a prime and let $a\in\Z$ with $p\nmid a$.
Then
\begin{equation}\label{tancot}
\begin{aligned}
&\sum_{k=1}^{(p-1)/2}\f1{\cot\pi\f{ak^2}p-1}=\sum_{k=1}^{(p-1)/2}\f1{1-\tan\pi\f{ak^2}p}-\f{p-1}2
\\=&\f p4\l(\l(\f{-1}p\r)-1\r)
+\l(\f{-2a}p\r)\f{\sqrt p}2\sum_{k=1}^{(p-1)/2}(-1)^k\l(\f{k}p\r).
\end{aligned}\end{equation}
\end{theorem}
\begin{remark}\label{Rem1.5} For any prime $p\eq1\pmod4$, we have
$$\sum_{k=1}^{(p-1)/2}\l(\f kp\r)=\f12\sum_{k=1}^{(p-1)/2}\l(\l(\f kp\r)+\l(\f{p-k}p\r)\r)=\sum_{r=1}^{p-1}\l(\f rp\r)=0$$
and hence
$$\sum_{k=1}^{(p-1)/2}(-1)^k\l(\f kp\r)
=\sum_{k=1}^{(p-1)/2}(1+(-1)^k)\l(\f kp\r)=2\sum_{j=1}^{\lfloor p/4\rfloor}\l(\f{2j}p\r)
=\l(\f 2p\r)h(-p)$$
(with $h(-p)$ the class number of the field $\Q(\sqrt{-p})$) since $\f{h(-p)}2=\sum_{0<k<p/4}(\f kp)$ (cf. \cite[(1.2)]{WC}), therefore we have the new class number formula:
\begin{equation}\label{h(-p)}
h(-p)=\f2{\sqrt p}\sum_{k=1}^{(p-1)/2}\f1{\cot\pi \f{k^2}p-1}.
\end{equation}
\end{remark}

\begin{theorem}\label{Th1.4} Let $p$ be an odd prime, and let $a\in\Z$ with $p\nmid a$.
Let $\ve_p$ and $h(p)$ be the fundamental unit and the class number of the real quadratic field
$\Q(\sqrt p)$.

{\rm (i)} If $p\eq1\pmod 8$, then
\begin{align}\label{tan1}
\prod_{k=1}^{(p-1)/2}\l(1+\tan\pi\f{ak^2}p\r)=&(-1)^{|\{1\ls k<\f p4:\ (\f kp)=1\}|}2^{(p-1)/4},
\\\label{cot1}
\prod_{k=1}^{(p-1)/2}\l(1+\cot\pi\f{ak^2}p\r)=&(-1)^{|\{1\ls k<\f p4:\ (\f kp)=1\}|}\f{2^{(p-1)/4}}{\sqrt p}\ve_p^{(\f ap)h(p)}.
\end{align}
If $p\eq5\pmod 8$, then
\begin{align}\label{tan5}
\prod_{k=1}^{(p-1)/2}\l(1+\tan\pi\f{ak^2}p\r)=&(-1)^{|\{1\ls k<\f p4:\ (\f kp)=-1\}|}2^{(p-1)/4}
\l(\f ap\r)\ve_p^{-3(\f ap)h(p)},
\\\label{cot5}
\prod_{k=1}^{(p-1)/2}\l(1+\cot\pi\f{ak^2}p\r)=&(-1)^{|\{1\ls k<\f p4:\ (\f kp)=1\}|}\l(\f ap\r)\f{2^{(p-1)/4}}{\sqrt p}.
\end{align}

{\rm (ii)} Suppose that $p\eq3\pmod 4$ and write $\varepsilon_p^{h(p)}=a_p+b_p\sqrt p$ with $a_p$ and $b_p$ positive integers. Set
$$s_p=\sqrt{a_p+(-1)^{(p+1)/4}}\ \ \ \text{and}\ \ \ t_p=\frac{b_p}{s_p}.$$
Then
\begin{equation}\label{tan43}
\prod_{k=1}^{(p-1)/2}\l(1+\tan\pi\f{ak^2}p\r)
=(-1)^{\da_{p,3}+\lfloor\f{p+1}8\rfloor+\f{h(-p)+1}2\cdot\f{p+1}4}2^{\f{p-3}4}\l(s_p+\l(\f ap\r)t_p\sqrt p\r),
\end{equation}
where the Kronecker symbol $\da_{p,3}$ takes $1$ or $0$ according as $p=3$ or not. Also,
\begin{equation}\label{cot43}
\prod_{k=1}^{(p-1)/2}\l(1+\cot\pi\f{ak^2}p\r)
=(-1)^{\lfloor\f{p-3}8\rfloor+\f{h(-p)-1}2\cdot\f{p-3}4}2^{\f{p-3}4}\l(t_p+\l(\f ap\r)\f{s_p}{\sqrt p}\r).
\end{equation}
\end{theorem}

Let $p$ be an odd prime, and define the polynomials
\begin{equation}\label{S-poly}G_p(x):=\prod_{k=1\atop (\f kp)=1}^{(p-1)/2}(x-e^{2\pi ik/p})=\prod_{k=1}^{(p-1)/2}(x-e^{2\pi ik^2/p})
\end{equation}
and
\begin{equation}\label{-S-poly}G_p^-(x):=\prod_{k=1\atop (\f kp)=-1}^{p-1}(x-e^{2\pi ik/p}).
\end{equation}
Clearly
$$G_p(x)G_p^-(x)=\prod_{r=1}^{p-1}(x-e^{2\pi ir/p})=\f{x^p-1}{x-1}.$$
In the case $p\eq3\pmod4$,  Dirichlet (cf. \cite[pp.\,370-371]{D99}) realized that
$(i-(\f 2p))G_p(i)\in\Z[\sqrt p]$, and Williams \cite[Lemma 3]{W82} determined
the exact value of $G_p(\pm i)$ which will be used in our proof of Theorem 1.4(ii). In the case $p\eq1\pmod8$, P. Kaplan and K. S. Williams \cite[(2.13)]{KW} proved that
$$G_p^-(e^{2\pi i/8})=(-1)^{(p-1)/8}i^{-(h(-p)+h(-2p))/4}\ve_{2p}^{h(2p)/8},$$
where $h(-p),h(-2p),h(2p)$ are class numbers of the quadratic fields
$\Q(\sqrt{-p}),\Q(\sqrt{-2p}),\Q(\sqrt{2p})$ respectively, and $\ve_{2p}$ is the fundamental unit
of the real quadratic field $\Q(\sqrt{2p})$.
In the case $p\eq5\pmod8$, the proof of the main theorem of Williams \cite{W81} involves $G_p^-(e^{2\pi i/8})$, but he did not present an explicit closed formula for it.
To prove Theorem 1.4(i), we will present in Section 4 a closed formula for $G_p(i)$ in the case $p\eq1\pmod4$.

Our next theorem determines the value of $G_p(\omega)$ for each prime $p\eq1\pmod4$, where
$$\omega=e^{2\pi i/3}=\f{-1+\sqrt3\,i}2.$$

\begin{theorem}\label{Th1.5} Let $p$ be a prime with $p\eq1\pmod 4$. Then
\begin{equation}\label{omega}(-1)^{|\{1\ls k\ls\lfloor\f{p+1}3\rfloor:\ (\f kp)=-1\}|}G_p(\omega)
=\begin{cases}1&\t{if}\ p\eq1\pmod{12},
\\\omega\ve_p^{h(p)}&\t{if}\ p\eq5\pmod{12}.
\end{cases}
\end{equation}
Also,
\begin{equation}\label{-omega41}G_p(-\omega)=\begin{cases}1&\t{if}\ p\eq1\pmod {12},
\\-\omega\ve_p^{-2h(p)}&\t{if}\ p\eq5\pmod{24},
\\\omega&\t{if}\ p\eq17\pmod{24}.
\end{cases}\end{equation}
\end{theorem}

We will show Theorem 1.1 and Corollary 1.1 in the next section, and prove Theorem 1.2 and Corollary 1.4 in Section 3. Our proof of Theorem 1.2 utilizes
the functional equation \eqref{csc}. Theorems 1.3-1.5 will be proved in Section 4.
We will pose some conjectures for further research in Section 5.

\section{Proofs of Theorem 1.1 and Corollary 1.1}
\setcounter{lemma}{0}
\setcounter{theorem}{0}
\setcounter{corollary}{0}
\setcounter{remark}{0}
\setcounter{equation}{0}

\begin{lemma}\label{Lem2.1} Let $n$ be any positive odd integer. Then
\begin{equation}\label{cot}\prod_{r=0}^{n-1}\l(1+\cot\pi\f{x+r}n\r)=\l(\f 2n\r)2^{(n-1)/2}\l(1+\l(\f{-1}n\r)\cot \pi x\r)
\end{equation}
for all $x\in\C\sm\Z$, and
\begin{equation}\label{tan}\prod_{r=0}^{n-1}\l(1+\tan\pi\f{x+r}n\r)=\l(\f 2n\r)2^{(n-1)/2}\l(1+\l(\f{-1}n\r)\tan \pi x\r)
\end{equation}
for all $x\in\C$ with $x-1/2\not\in\Z$.
\end{lemma}
\Proof. Let $x\in\C\sm\Z$. For each $r=0,\ldots,n-1$, by Euler's formula $e^{iz}=\cos z+i\sin z$ we have
\begin{align*}1+\cot\pi\f{x+r}n
=&1+\f{(e^{i\pi(x+r)/n}+e^{-i\pi(x+r)/n})/2}{(e^{i\pi(x+r)/n}-e^{-i\pi(x+r)/n})/(2i)}
\\=&1+i\f{e^{2\pi i(x+r)/n}+1}{e^{2\pi i(x+r)/n}-1}=1+i\l(1+\f2{e^{2\pi i(x+r)/n}-1}\r)
\\=&(1+i)\l(1+\f{1+i}{e^{2\pi i(x+r)/n}-1}\r)=(1+i)\f{-i-e^{2\pi i(x+r)/n}}{1-e^{2\pi i(x+r)/n}}.
\end{align*}
Thus
\begin{align*}\prod_{r=0}^{n-1}\l(1+\cot\pi\f{x+r}n\r)
=&(1+i)^n\f{\prod_{r=0}^{n-1}(y-e^{2\pi i(x+r)/n})\big|_{y=-i}}{\prod_{r=0}^{n-1}(z-e^{2\pi i(x+r)/n})\big|_{z=1}}
\\=&(1+i)\l((1+i)^2\r)^{(n-1)/2}\f{(y^n-e^{2\pi ix})|_{y=-i}}{(z^n-e^{2\pi ix})|_{z=1}}
\\=&(1+i)(2i)^{(n-1)/2}\f{e^{2\pi ix}+i(-1)^{(n-1)/2}}{e^{2\pi ix}-1}.
\end{align*}
On the other hand,
\begin{align*}1+\l(\f{-1}n\r)\cot\pi x=&1+(-1)^{(n-1)/2}i\f{e^{i\pi x}+e^{-i \pi x}}
{e^{i\pi x}-e^{-i\pi x}}
\\=&1+(-1)^{(n-1)/2}i\f{e^{2\pi ix}+1}{e^{2\pi ix}-1}
\\=&(1+(-1)^{(n-1)/2}i)
\f{e^{2\pi ix}+i(-1)^{(n-1)/2}}{e^{2\pi ix}-1}.
\end{align*}
Therefore
$$\prod_{r=0}^{n-1}\l(1+\cot\pi\f{x+r}n\r)
=\f{(1+i)i^{(n-1)/2}}{1+(-1)^{(n-1)/2}i}2^{(n-1)/2}\l(1+\l(\f{-1}n\r)\cot\pi x\r).$$
Since
$$\f{(1+i)i^{(n-1)/2}}{1+(-1)^{(n-1)/2}i}=(-1)^{(n^2-1)/8}=\l(\f 2n\r),$$
we obtain \eqref{cot} from the above.

Now let $x\in\C$ with $x-1/2\not\in\Z$. Then $x'=n/2-x\not\in\Z$.
Applying \eqref{cot} with $x$ replaced by $x'$, we get that
$$\prod_{r=0}^{n-1}\l(1+\cot\pi\f{x'+r-n}n\r)=\l(\f 2n\r)2^{(n-1)/2}\l(1+\l(\f{-1}n\r)\cot\l(n\f{\pi}2-\pi x\r)\r),$$
i.e.,
$$\prod_{r=0}^{n-1}\l(1+\tan\pi\f{x+(n-r)}n\r)=\l(\f 2n\r)2^{(n-1)/2}\l(1+\l(\f{-1}n\r)\tan\pi x\r).$$
Therefore \eqref{tan} holds. \qed

\medskip
\noindent{\it Proof of Theorem 1.1}. Observe that
$$\f{\f d{dz}(1+\tan z)}{1+\tan z}=\f{\sec^2 z}{1+\tan z}=\f1{\cos^2z+\sin z\cos z}=\f 2{1+\cos 2z+\sin 2z}.$$
By taking the logarithmic derivative, we obtain from \eqref{tan} the equality
$$\sum_{r=0}^{n-1}\f{2\pi/n}{1+\cos2\pi(x+r)/n+\sin2\pi(x+r)/n}=\f{(\f{-1}n)2\pi}{1+\cos2\pi x+(\f{-1}n)\sin2\pi x},$$
provided that $x+1/2,x+(-1)^{(n-1)/2}/4\not\in\Z$. (Note that
$1+(\f{-1}n)\tan\pi x=0$ if and only if $x+(\f{-1}n)\f14\in\Z$.)
This proves the first assertion in Theorem 1.1.

Now let $x\in\C$ with $x\not\in\Z$ and $x+(-1)^{(n-1)/2}/4\not\in\Z$.
Set $x'=n/2-x$. Then $x'+1/2\not\in\Z$ and $x'+(-1)^{(n-1)/2}/4\not\in\Z$. By \eqref{sincos} with $x$ replaced by $x'$, we have
\begin{align*}&\sum_{r=0}^{n-1}\f1{1+\sin(\pi-2\pi(x+n-r)/n)+\cos(\pi-2\pi(x+n-r)/n)}
\\=&\f{(\f{-1}n)n}{1+(\f{-1}n)\sin(n\pi-2\pi x)+\cos(n\pi-2\pi x)},
\end{align*}
i.e.,
$$\sum_{r=0}^{n-1}\f1{1+\sin2\pi\f{x+(n-r)}n-\cos2\pi\f{x+(n-r)}n}=\f{(\f{-1}n)n}{1+(\f{-1}n)\sin2\pi x-\cos2\pi x}.$$
Therefore \eqref{-sincos} holds. This proves the second assertion in Theorem 1.1. \qed

\medskip
\noindent
{\it Proof of Corollary 1.1}. Let $x\in\C$ with $2x\not\in\Z$. We want to show \eqref{csc}.

We first assume that $4x$ is not an odd integer. Then
both \eqref{sincos} and \eqref{-sincos} hold, and hence
\begin{align*}&\sum_{r=0}^{n-1}\l(\f1{1+\sin2\pi\f{x+r}n+\cos2\pi\f{x+r}n}+\f1{1+\sin2\pi\f{x+r}n-\cos2\pi\f{x+r}n}\r)
\\=&\l(\f{-1}n\r)n\l(\f1{1+(\f{-1}n)\sin 2\pi x+\cos 2\pi x}+\f1{1+(\f{-1}n)\sin 2\pi x-\cos 2\pi x}\r),
\end{align*}
i.e.,
\begin{align*}&\sum_{r=0}^{n-1}\f{2(1+\sin2\pi\f{x+r}n)}{(1+\sin2\pi\f{x+r}n)^2-(1-\sin^22\pi\f{x+r}n)}
\\=&\l(\f{-1}n\r)n\f{2(1+(\f{-1}n)\sin2\pi x)}{(1+(\f{-1}n)\sin2\pi x)^2-(1-\sin^2 2\pi x)}.
\end{align*}
Therefore \eqref{csc} follows.

When $4x$ is an odd integer, by the above we have
$$\sum_{r=0}^{n-1}\csc2\pi\f{x+r}n=\lim_{t\to x\atop 4t\not\in\Z}\sum_{r=0}^{n-1}\csc2\pi\f{t+r}n
=\lim_{t\to x\atop 4t\not\in\Z}n\csc2\pi t=n\csc2\pi x.$$

Now we turn to show \eqref{sec} for any complex number $x$ with $4x$ not an odd integer.
For $x'=n/4-x$, we have $2x'=n/2-2x\not\in\Z$. Applying \eqref{csc} with $x$ replaced by $x'=n/4-x$,
we find that
$$\f1n\sum_{r=0}^{n-1}\csc\l(\f{\pi}2-2\pi\f{x+n-r}n\r)=\csc\l(n\f{\pi}2-2\pi x\r)$$
and hence
$$\f1n\sum_{r=0}^{n-1}\sec2\pi\f{x+(n-r)}n=(-1)^{(n-1)/2}\sec2\pi x,$$
which is equivalent to \eqref{sec}. This concludes the proof. \qed

\section{Proofs of Theorem 1.2 and Corollary 1.4}
\setcounter{lemma}{0}
\setcounter{theorem}{0}
\setcounter{corollary}{0}
\setcounter{remark}{0}
\setcounter{equation}{0}

\medskip
\noindent{\it Proof of Theorem 1.2}. Let $x,y\in\C$ with $x+y\not\in\Z$
and $x-y-1/2\not\in\Z$.
For any $j,k=0,\ldots,n-1$, clearly $2(x-y+j-k)/n$ cannot be an odd integer and thus
\begin{align*}&2i\sin 2\pi\f{x+j}n+2i\sin2\pi\f{y+k}n
\\=&e^{2\pi i(x+j)/n}-e^{-2\pi i(x+j)/n}+e^{2\pi i(y+k)/n}-e^{-2\pi i(y+k)/n}
\\=&(e^{2\pi i(x+j)/n}+e^{2\pi i(y+k)/n})(1-e^{-2\pi i(x+y+j+k)/n})\not=0.
\end{align*}
As $n$ is odd, $\{2r:\ r=0,\ldots,n-1\}$ is a complete system of residues modulo $n$.
Let $L$ denote the left-hand side of \eqref{cos}. By the above,
\begin{align*}L=&\sum_{j,k=0}^{n-1}\f{2i}{(e^{2\pi i(x+j)/n}+e^{2\pi i(y+k)/n})(1-e^{-2\pi i(x+y+j+k)/n})}
\\=&\sum_{r=0}^{n-1}\f{i}{1-e^{-2\pi i(x+y+2r)/n}}\sum_{k=0}^{n-1}\f2{e^{2\pi i(y+k)/n}+e^{2\pi i(x+2r-k)/n}}
\\=&\sum_{r=0}^{n-1}\f{ie^{-2\pi iy/n}}{1-e^{-2\pi i(x+y+2r)/n}}\sum_{k=0}^{n-1}\f2{e^{2\pi ik/n}+e^{2\pi i(x-y+2r-k))/n}}
\\=&\sum_{r=0}^{n-1}\f{ie^{-2\pi iy/n}}{1-e^{-2\pi i(x+y+2r)/n}}\sigma_r,
\end{align*}
where
$$\sigma_r:=\sum_{k=0}^{n-1}\l(\f1{e^{2\pi ik/n}+ie^{2\pi i((x-y)/2+r)/n}}+\f1{e^{2\pi ik/n}-ie^{2\pi i((x-y)/2+r)/n}}\r).$$

As $\prod_{k=0}^{n-1}(z-e^{2\pi ik/n})=z^n-1,$
by taking the logarithmic derivative we get
$$\sum_{k=0}^{n-1}\f1{z-e^{2\pi ik/n}}=\f{nz^{n-1}}{z^n-1},\ \t{i.e.},\ \sum_{k=0}^{n-1}\f1{e^{2\pi ik/n}-z}=\f{nz^{n-1}}{1-z^n}.$$
Hence, for each $r=0,\ldots,n-1$ we have
\begin{align*}\sigma_r=&\f{n(-ie^{2\pi i((x-y)/2+r)/n})^{n-1}}{1-(-ie^{2\pi i((x-y)/2+r)/n})^n}+\f{n(ie^{2\pi i((x-y)/2+r)/n})^{n-1}}{1-(ie^{2\pi i((x-y)/2+r)/n})^n}
\\=&n(-1)^{(n-1)/2}e^{-2\pi ir/n}e^{i\pi(x-y)(n-1)/n}\l(\f1{1+i^ne^{i\pi(x-y)}}
+\f1{1-i^ne^{i\pi(x-y)}}\r)
\\=&n(-1)^{(n-1)/2}e^{-2\pi ir/n}e^{i\pi(x-y)(n-1)/n}\f2{1+e^{2\pi i(x-y)}}
\\=&n(-1)^{(n-1)/2}e^{-2\pi ir/n}e^{-i\pi(x-y)/n}\f2{e^{-i\pi (x-y)}+e^{i\pi (x-y)}}
\end{align*}

In view of the above, we see that
\begin{align*}L=&\sum_{r=0}^{n-1}\f{ie^{-2\pi iy/n}}{1-e^{-2\pi i(x+y+2r)/n}}n(-1)^{(n-1)/2}e^{-2\pi ir/n}\f{2e^{-i\pi(x-y)/n}}{e^{-i\pi (x-y)}+e^{i\pi (x-y)}}
\\=&\f{(-1)^{(n-1)/2}2n}{e^{i\pi (x-y)}+e^{-i\pi (x-y)}}\sum_{r=0}^{n-1}\f{ie^{-2\pi i((x+y)/2+r)/n}}{1-e^{-2\pi i(x+y+2r)/n}}
\\=&\f{(-1)^{(n-1)/2}n}{e^{i\pi (x-y)}+e^{-i\pi (x-y)}}\sum_{r=0}^{n-1}\f {2i}{e^{2\pi i((x+y)/2+r)/n}-e^{-2\pi i((x+y)/2+r)/n}}
\\=&\f{(-1)^{(n-1)/2}n}{e^{i\pi (x-y)}+e^{-i\pi (x-y)}}\sum_{r=0}^{n-1}\csc2\pi\f{(x+y)/2+r}n.
\end{align*}
Combining this with \eqref{csc}, we obtain
\begin{align*}L=&\f{(-1)^{(n-1)/2}n}{e^{i\pi (x-y)}+e^{-i\pi (x-y)}}\times n\f{2i}{e^{i\pi (x+y)}-e^{-i\pi (x+y)}}
\\=&\f{(\f{-1}n)n^2\times 2i}{e^{2\pi ix}-e^{-2\pi ix}+e^{2\pi iy}-e^{-2\pi iy}}
=\l(\f{-1}n\r)\f{n^2}{\sin2\pi x+\sin2\pi y}.
\end{align*}
So \eqref{sin} holds.

Now let $x,y\in\C$ with $x\pm y-1/2\not\in\Z$. For $x'=n/4-x$ and $y'=n/4-y$, we have
$x'+y'=n/2-(x+y)\not\in\Z$ and $x'-y'-1/2=y-x-1/2\not\in\Z$. Thus, by \eqref{sin} with $x$ and $y$ replaced by $x'$ and $y'$ respectively, we get
$$\sum_{j,k=0}^{n-1}\f1{\sin 2\pi\f{x'+j-n}n+\sin2\pi\f{y'+k-n}n}
=\l(\f{-1}n\r)\f{n^2}{\sin(n\f{\pi}2-2\pi x)+\sin(n\f{\pi}2-2\pi y)},$$
i.e.,
$$\sum_{j,k=0}^{n-1}\f1{\cos 2\pi(x+(n-j))/n+\cos2\pi(y+(n-k))/n}=\f{n^2}{\cos2\pi x+\cos2\pi y}.$$
Therefore \eqref{cos} holds.

Finally, we let $x,y\in\C$ with $x\pm y+(-1)^{(n-1)/2}/4\not\in\Z$.
Set $x'=n/4-x$. Then $x'\pm y-1/2=n/4-1/2-x\pm y\not\in\Z$. Applying \eqref{cos} with $x$ replaced by
$x'$, we get
$$\sum_{j,k=0}^{n-1}\f1{\cos (\f{\pi}2-2\pi\f{x+(n-j)}n)+\cos2\pi\f{y+k}n}=\f{n^2}{\cos(n\f{\pi}2-2\pi x)+\cos2\pi y},$$
which is equivalent to \eqref{mix}.

The proof of Theorem \ref{Th1.2} is now complete. \qed

\medskip
\noindent{\it Proof of Corollary 1.4}. As $(\f{-1}p)=-1$, we see that
$$\l\{1^2,\ldots,\l(\f{p-1}2\r)^2,-1^2,\ldots,-\l(\f{p-1}2\r)^2\r\}$$
is a reduced system of residue classes modulo $p$. Thus
\begin{align*}&\sum_{s,t=1}^{p-1}\f1{\cos2\pi s/p+\cos2\pi t/p}
\\=&\sum_{s=1}^{p-1}\sum_{k=1}^{(p-1)/2}\l(\f1{\cos2\pi s/p+\cos2\pi k^2/p}+\f1{\cos2\pi s/p+\cos2\pi (-k^2)/p}\r)
\\=&2\sum_{s=1}^{p-1}\sum_{k=1}^{(p-1)/2}\f1{\cos2\pi s/p+\cos2\pi k^2/p}=4\sum_{j,k=1}^{(p-1)/2}
\f1{\cos2\pi j^2/p+\cos2\pi k^2/p}
\\=&4\sum_{j=1}^{(p-1)/2}\f1{2\cos2\pi j^2/p}+8\sum_{1\ls j<k\ls(p-1)/2}\f1{\cos2\pi j^2/p+\cos2\pi k^2/p}
\\=&\sum_{r=1}^{p-1}\sec2\pi \f rp+8\sum_{1\ls j<k\ls(p-1)/2}\f1{\cos2\pi j^2/p+\cos2\pi k^2/p}
\end{align*}
and hence
\begin{align*}&\sum_{s,t=0}^{p-1}\f1{\cos2\pi s/p+\cos2\pi t/p}-8\sum_{1\ls j<k\ls(p-1)/2}\f1{\cos2\pi j^2/p+\cos2\pi k^2/p}
\\=&\sum_{s=0}^{p-1}\f1{\cos2\pi s/p+\cos0}+\sum_{t=0}^{p-1}\f1{\cos0+\cos2\pi t/p}-\f1{2\cos0}
+\sum_{r=1}^{p-1}\sec2\pi\f rp
\\=&\sum_{r=0}^{p-1}\f2{1+\cos2\pi r/p}-\f12+\sum_{r=1}^{p-1}\sec2\pi \f rp
=\sum_{r=0}^{p-1}\sec^2\pi\f rp+\sum_{r=0}^{p-1}\sec2\pi \f rp-\f 32.
\end{align*}
With the help of \eqref{cos0}, \eqref{secant} and \eqref{sec0}, we finally obtain
\begin{align*}&8\sum_{1\ls j<k\ls(p-1)/2}\f1{\cos2\pi j^2/p+\cos2\pi k^2/p}
\\=&\f{p^2}2-p^2-\l(\f{-1}p\r)p+\f 32=-\f{(p+1)(p-3)}2
\end{align*}
and hence the desired identity \eqref{cosp} follows. \qed
\medskip

\section{Proofs of Theorems 1.3-1.5}
\setcounter{lemma}{0}
\setcounter{theorem}{0}
\setcounter{corollary}{0}
\setcounter{remark}{0}
\setcounter{equation}{0}

\medskip
\noindent
{\it Proof of Theorem 1.3}. Clearly,
\begin{align*}&\sum_{k=1}^{(p-1)/2}\f1{1-\tan \pi ak^2/p}+\sum_{k=1}^{(p-1)/2}\f1{1-\cot\pi ak^2/p}
\\=&\sum_{k=1}^{(p-1)/2}\l(\f{\cos\pi ak^2/p}{\cos\pi ak^2/p-\sin\pi ak^2/p}-\f{\sin\pi ak^2/p}{\cos\pi ak^2/p-\sin\pi ak^2/p}\r)
\\=&\sum_{k=1}^{(p-1)/2}1=\f{p-1}2.
\end{align*}
So the first equality in \eqref{tancot} holds.

For any $x\in\C$ with $2x$ not an odd integer, we have
$$\cot\pi x-1=i\f{e^{2\pi ix}+1}{e^{2\pi ix}-1}-1=(i-1)\f{e^{2\pi ix}-i}{e^{2\pi ix}-1}.$$
Thus
\begin{align*}\sum_{k=1}^{(p-1)/2}\f1{\cot\pi\f{ak^2}p-1}=&\f1{i-1}\sum_{k=1}^{(p-1)/2}\l(1+\f{i-1}{e^{2\pi iak^2/p}-i}\r)
\\=&\f1{i-1}\cdot\f{p-1}2+\f1{1-i^p}\sum_{k=1}^{(p-1)/2}\f{(e^{2\pi iak^2/p})^p-i^p}{e^{2\pi i ak^2/p}-i}.
\end{align*}
Observe that
\begin{align*}&\sum_{k=1}^{(p-1)/2}\f{(e^{2\pi iak^2/p})^p-i^p}{e^{2\pi i ak^2/p}-i}
=\sum_{k=1}^{(p-1)/2}\sum_{j=0}^{p-1}(e^{2\pi iak^2/p})^ji^{p-1-j}
\\=&\f{p-1}2i^{p-1}+\sum_{j=1}^{p-1}\f{i^{p-1-j}}2\sum_{k=1}^{(p-1)/2}\l(e^{2\pi ijak^2/p}+e^{2\pi ija(p-k)^2/p}\r)
\\=&\f{p-1}2i^{p-1}+\sum_{j=1}^{p-1}\f{i^{p-1-j}}2\(\sum_{x=0}^{p-1}e^{2\pi iajx^2/p}-1\)
\\=&\f p2i^{p-1}-\f12\sum_{j=0}^{p-1}i^{p-1-j}+\f12\sum_{j=1}^{p-1}i^{p-1-j}\l(\f {aj}p\r)\sqrt{(-1)^{(p-1)/2}p}
\end{align*}
with the help of the known evaluations of quadratic Gauss sums.
Therefore
\begin{align*}\sum_{k=1}^{(p-1)/2}\f1{\cot\pi\f{ak^2}p-1}
=&\f1{i-1}\cdot\f{p-1}2+\f1{1-i^p}\l(\f p2 i^{p-1}-\f12\cdot\f{1-i^p}{1-i}\r)
\\&+\f{\sqrt{(-1)^{(p-1)/2}p}}{2(1-i^p)}\sum_{j=1}^{p-1}i^{p-1-j}\l(\f {aj}p\r)
\\=&\f p2\l(\f1{i-1}+\f{i^{p-1}}{1-i^p}\r)
+\l(\f ap\r)\f{\sqrt{(-1)^{(p-1)/2}p}}{2(1-i^p)}
\\&\times\sum_{k=1}^{(p-1)/2}\l(i^{p-1-2k}\l(\f{2k}p\r)+i^{p-1-(p-2k)}\l(\f{p-2k}p\r)\r)
\end{align*}
and hence
\begin{align*}&\sum_{k=1}^{(p-1)/2}\f1{\cot\pi\f{ak^2}p-1}+\f p2\cdot\f{1-i^{p-1}}{(1-i)(1-i^p)}
\\=&\l(\f ap\r)\f{\sqrt{(-1)^{(p-1)/2}p}}{2(1-i^p)}(1-i)\sum_{k=1}^{(p-1)/2}(-1)^k\l(\f{-2k}p\r).
\end{align*}
It follows that
$$\sum_{k=1}^{(p-1)/2}\f1{\cot\pi\f{ak^2}p-1}=\f p4\l(\l(\f{-1}p\r)-1\r)
+\l(\f{-2a}p\r)\f{\sqrt p}2\sum_{k=1}^{(p-1)/2}(-1)^k\l(\f{k}p\r).$$

In view of the above, we have proved \eqref{tancot} fully. \qed

\begin{lemma}\label{Lem4.1} {\rm (Williams and Currie \cite[(1.4)]{WC})}
Let $p\eq1\pmod4$ be a prime. Then
\begin{equation}\label{wc} (-1)^{|\{1\ls k< \f p4:\ (\f kp)=-1\}|}2^{(p-1)/4}\eq\begin{cases}1\pmod p&\t{if}\ p\eq1\pmod8,
\\\f{p-1}2!\pmod p&\t{if}\ p\eq5\pmod 8.\end{cases}\end{equation}
\end{lemma}

\begin{lemma}\label{Lem4.2} Let $p\eq1\ (\mo\ 4)$ be a prime.

{\rm (i) (\cite[(1.12) and (1.17)]{S19b})} For any integer $a\not\eq0\pmod p$, we have
\begin{equation}\label{p14}\prod_{k=1}^{(p-1)/2}(1-e^{2\pi iak^2/p})=\sqrt p\,\ve_p^{-(\f ap)h(p)}
\end{equation}
and
\begin{equation}\label{cos14}2^{(p-1)/2}\prod_{k=1}^{(p-1)/2}\cos\pi \f{ak^2}p=
(-1)^{a(p-1)/4}\ve_p^{(1-(\f 2p))(\f ap)h(p)}
\end{equation}

{\rm (ii) (\cite[Corollary 1.1]{S19a})}
 Write $\ve_p^{h(p)}=a_p+b_p\sqrt p$ with $2a_p,2b_p\in\Z$.
Then
\begin{equation}\label{re}a_p\eq-\f{p-1}2!\pmod p.
\end{equation}
\end{lemma}

\begin{lemma}\label{Lem-4.3} Let $p>3$ be a prime with $p\eq3\pmod4$. For any integer $a\not\eq0\pmod p$,
we have
\begin{equation}\label{4k+3}\prod_{k=1}^{(p-1)/2}(1-e^{2\pi iak^2/p})=(-1)^{(h(-p)+1)/2}\l(\f ap\r)\sqrt p\,i.
\end{equation}
\end{lemma}
\begin{remark} This follows from Dirichlet's class number formula, see Williams \cite[Lemma 3]{W81}
and Sun \cite[Theorem 1.3(i)]{S19b}.
\end{remark}

\begin{lemma} \label{Lem-S_p(i)} Let $p>3$ be a prime.

{\rm (i)} Let $a$ be any integer not divisible by $p$. If $p\eq1\pmod 8$, then
\begin{equation}\label{root1}\prod_{k=1}^{(p-1)/2}(i-e^{2\pi iak^2/p})=(-1)^{\f{p-1}8+|\{1\le k<\f p4:\ (\f kp)=1\}|}.
\end{equation}
If $p\eq5\pmod 8$, then
\begin{equation}\label{root5}\prod_{k=1}^{(p-1)/2}(i-e^{2\pi iak^2/p})=i(-1)^{\f{p-5}8+|\{1\le k<\f p4:\ (\f kp)=1\}|}
\l(\f ap\r)\ve_p^{-(\f ap)h(p)}.
\end{equation}

{\rm (ii)} When $p\eq3\pmod4$, we have
\begin{equation}\label{S(i)}\left(i-(-1)^{(p+1)/4}\right)G_p(i)=(-1)^{\frac{h(-p)+1}2\cdot\frac{p+1}4}(s_p-t_p\sqrt p)\end{equation}
with $s_p$ and $t_p$ given by Theorem 1.4.
\end{lemma}
\Proof. (i) Suppose that $p\eq1\pmod4$.
Let $c:=\prod_{k=1}^{(p-1)/2}(i-e^{2\pi iak^2/p}).$
In the ring of algebraic $p$-adic integers, we have the congruence
\begin{equation}\label{cp}c^p\eq\prod_{k=1}^{(p-1)/2}(i^p-1)=(i-1)^{(p-1)/2}
=(-2i)^{(p-1)/4}\pmod p.
\end{equation}

As $(\f{-1}p)=1$, we have
\begin{align*}c^2=&\prod_{k=1}^{(p-1)/2}\l(i-e^{2\pi iak^2/p}\r)\l(i-e^{-2\pi iak^2/p}\r)
=\prod_{k=1}^{(p-1)/2}\l(-ie^{2\pi iak^2/p}-ie^{-2\pi iak^2/p}\r)
\\=&(2i)^{(p-1)/2}\prod_{k=1}^{(p-1)/2}\cos\pi\f{2ak^2}p=(-1)^{(p-1)/4}\ve_p^{(1-(\f 2p))(\f {2a}p)h(p)}
\end{align*}
with the aid of \eqref{cos14}, and hence
\begin{equation}\label{cc}c=\da i^{(p-1)/4}\ve_p^{((\f 2p)-1)(\f ap)h(p)/2}
\end{equation}
for some $\da\in\{\pm1\}$. Note that $i^p=i$. Thus
$$c^p= \da i^{(p-1)/4}\ve_p^{((\f 2p)-1)(\f ap)ph(p)/2}.$$
Combining this with \eqref{cp} we get
$$(-2)^{(p-1)/4}\eq \da\ve_p^{((\f 2p)-1)(\f ap)ph(p)/2}\pmod p.$$
If $p\eq1\pmod 8$, then
$$\da\eq 2^{(p-1)/4}\eq(-1)^{|\{1\ls k<\f p4:\ (\f kp)=-1\}|}=(-1)^{|\{1\ls k<\f p4:\ (\f kp)=1\}|}\pmod p$$ with the aid of \eqref{wc}.
When $p\eq5\pmod 8$, we write $\ve_p^{h(p)}=a_p+b_p\sqrt p$ with $2a_p,2b_p\in\Z$, and observe that
\begin{align*}
(-1)^{|\{1\ls k<\f p4:\ (\f kp)=1\}|}\da\f{p-1}2!
\eq&-\da2^{(p-1)/4}\eq\ve_p^{-(\f ap)ph(p)}=(a_p+b_p\sqrt p)^{-(\f ap)p}
\\\eq&\l(a_p^p+b_p^pp^{(p-1)/2}\sqrt p\r)^{-(\f ap)}
\\\eq&\l(-\f{(p-1)!}2\r)^{-(\f ap)}\eq\l(\f ap\r)\f{p-1}2!\pmod p
\end{align*}
in light of \eqref{wc}, \eqref{re} and the simple fact that
$$\l(\f{p-1}2!\r)^2\eq\prod_{k=1}^{(p-1)/2}k(p-k)=(p-1)!\eq-1\pmod p$$
by Wilson's theorem. Therefore
$$\da=(-1)^{|\{1\ls k<p/4:\ (\f kp)=1\}|}\l(\f ap\r)^{(1-(\f 2p))/2}.$$
Combining this with \eqref{cc}, we obtain \eqref{root1} and \eqref{root5} in the cases $p\eq1\pmod8$
and $p\eq5\pmod 8$ respectively.

(ii) Now we handle the case $p\eq3\pmod4$. By Williams \cite[Lemma 3]{W81},
$$G_p(i)=-e^{2\pi i(\f 2p)/8}(-1)^{\f{h(-p)+1}2\cdot\f{p+1}4}\ve_p^{-h(p)/2}$$
and hence
$$\l(i-(-1)^{(p+1)/4}\r)G_p(i)=\l(\f 2p\r)(-1)^{\f{h(-p)+1}2\cdot\f{p+1}4}\sqrt{\f 2{a_p+b_p\sqrt p}}.$$
It is easy to see that
\begin{equation}\label{st} a_p^2-pb_p^2=1\ \t{and}\ \f{s_p^2-pt_p^2}2=\l(\f 2p\r)=G_p(i)G_p(-i)
\end{equation}
with the aid of \cite[(1.13)]{S19b}. Thus
\begin{align*}\sqrt{\f 2{a_p+b_p\sqrt p}}=&\sqrt{2(a_p-b_p\sqrt p)}
=\sqrt{2(s_p^2-\f{s_p^2-pt_p^2}2)-2s_pt_p\sqrt p}
\\=&\sqrt{(s_p-t_p\sqrt p)^2}=\l(\f 2p\r)(s_p-t_p\sqrt p)
\end{align*}
and hence \eqref{S(i)} holds.

In view of the above, we have completed the proof of Lemma \ref{Lem-S_p(i)}. \qed
\medskip

{\it Example} 4.1. For the prime $p=79$, we have $h(-p)=5$, $h(p)=3$ and $\varepsilon_p=80+9\sqrt p$. Note that
$$\varepsilon_p^{h(p)}=(80+9\sqrt{79})^3=2047760 + 230391\sqrt{79},$$
and
$$s_p=\sqrt{2047760+1}=1431\ \ \ \text{and}\ \ \ t_p= \frac{230391}{1431}=161.$$
Thus Lemma \ref{Lem-S_p(i)}(ii) for $p=79$ states that
$$(i-1)S_{79}(i)=1431-161\sqrt{79}.$$

\medskip
\noindent{\it Proof of Theorem 1.4}. From the proof of Lemma 2.1, for each $k=1,\ldots,(p-1)/2$
we have
$$1+\cot\pi\f{ak^2}p=(1+i)\f{-i-e^{2\pi iak^2/p}}{1-e^{2\pi iak^2/p}}.$$
Thus
\begin{equation}\label{cot-p}\prod_{k=1}^{(p-1)/2}\l(1+\cot\pi\f{ak^2}p\r)
=(1+i)^{(p-1)/2}\f{\prod_{k=1}^{(p-1)/2}\overline{i-e^{-2\pi iak^2/p}}}{\prod_{k=1}^{(p-1)/2}(1-e^{2\pi iak^2/p})},
\end{equation}
where $\bar z$ denotes the conjugate of $z\in\C$.

For each $k=1,\ldots,(p-1)/2$, we clearly have
\begin{align*}1+\tan\pi\f{ak^2}p=&1+\f{(e^{i\pi ak^2/p}-e^{-i\pi ak^2/p})/(2i)}
{(e^{i\pi ak^2/p}+e^{-i\pi ak^2/p})/2}
\\=&1+i\f{1-e^{2\pi iak^2/p}}{1+e^{2\pi iak^2/p}}=(1-i)\f{i+e^{2\pi iak^2/p}}{1+e^{2\pi iak^2/p}}.
\end{align*}
Thus
\begin{equation}\label{tan-p}
\prod_{k=1}^{(p-1)/2}\l(1+\tan\pi\f{ak^2}p\r)=(i-1)^{(p-1)/2}\f{\prod_{k=1}^{(p-1)/2}\overline{i-e^{-2\pi iak^2/p}}}{\prod_{k=1}^{(p-1)/2}(1+e^{2\pi iak^2/p})}.
\end{equation}

(i) Assume that $p\eq1\pmod4$. Combining \eqref{cot-p} and \eqref{p14}, we get
\begin{align*}\prod_{k=1}^{(p-1)/2}\l(1+\cot\pi\f{ak^2}p\r)
=&\f{(2i)^{(p-1)/4}}{\sqrt p\ve_p^{-(\f ap)h(p)}}\prod_{k=1}^{(p-1)/2}\overline{i-e^{-2\pi iak^2/p}}
\\=&\f{(2i)^{(p-1)/4}}{\sqrt p}\ve_p^{(\f ap)h(p)}\prod_{j=1}^{(p-1)/2}\overline{i-e^{2\pi iaj^2/p}}
\end{align*}
since $(\f{-1}p)=1$.
Therefore, \eqref{root1} implies \eqref{cot1} if $p\eq1\pmod8$,
and \eqref{root5} implies \eqref{cot5} if $p\eq5\pmod8$.

By Lemma 4.2(i),
\begin{align*}\prod_{k=1}^{(p-1)/2}\l(1+e^{2\pi iak^2/p}\r)
=&\f{\prod_{k=1}^{(p-1)/2}(1-e^{2\pi i2ak^2/p})}{\prod_{k=1}^{(p-1)/2}(1-e^{2\pi iak^2/p})}
\\=&\f{\sqrt p\,\ve_p^{-(\f{2a}p)h(p)}}{\sqrt p\,\ve_p^{-(\f ap)h(p)}}=\ve_p^{(1-(\f 2p))(\f ap)h(p)}.
\end{align*}
Combining this with \eqref{tan-p} and noting $(\f{-1}p)=1$, we obtain
\begin{align*}
\prod_{k=1}^{(p-1)/2}\l(1+\tan\pi\f{ak^2}p\r)
=(-2i)^{(p-1)/4}\ve_p^{((\f 2p)-1)(\f ap)h(p)}\prod_{k=1}^{(p-1)/2}\overline{i-e^{2\pi iak^2/p}}.
\end{align*}
Thus, \eqref{root1} implies \eqref{tan1} if $p\eq1\pmod8$,
and \eqref{root5} implies \eqref{tan5} if $p\eq5\pmod8$.

(ii) Now we handle the case $p\eq3\pmod4$. For $p=3$ we can easily verify \eqref{tan43} and \eqref{cot43}. Below we assume that $p>3$.

If $(\f ap)=-1$, then $(\f{-a}p)=1$ and hence
\begin{align*}\prod_{k=1}^{(p-1)/2}(i-e^{-2\pi iak^2/p})=&G_p(i)
=\f{(-1)^{\f{h(-p)+1}2\cdot\f{p+1}4}}{i+(-1)^{(p-3)/4}}(s_p-t_p\sqrt p)
\\=&\f{(-1)^{(p-3)/4}-i}2\times(-1)^{\f{h(-p)+1}2\cdot\f{p+1}4}(s_p-t_p\sqrt p)
\end{align*}
with the aid of \eqref{S(i)}.
Similarly, if $(\f ap)=1$, then $(\f{-a}p)=-1$ and hence
\begin{align*}\prod_{k=1}^{(p-1)/2}(i-e^{-2\pi iak^2/p})=&G_p^-(i)=\f{(i^p-1)/(i-1)}{G_p(i)}
\\=&\f{(i-(-1)^{(p+1)/4})(i^3-1)/(i-1)}
{(-1)^{\f{h(-p)+1}2\cdot\f{p+1}4}(s_p-t_p\sqrt p)}
\\=&\f{i(i+(-1)^{(p-3)/4})}{(-1)^{\f{h(-p)+1}2\cdot\f{p+1}4}}\times\f{s_p+t_p\sqrt p}{2(\f 2p)}
\\=&\f{(-1)^{(p-3)/4}-i}2\times(-1)^{\f{h(-p)+1}2\cdot\f{p+1}4}(s_p+t_p\sqrt p)
\end{align*}
with the aid of \eqref{st}.
Therefore
\begin{equation}\label{conjugate}\prod_{k=1}^{(p-1)/2}\overline{i-e^{-2\pi iak^2/p}}=\f{(-1)^{(p-3)/4}+i}2\times(-1)^{\f{h(-p)+1}2\cdot\f{p+1}4}\l(s_p+\l(\f ap\r)t_p\sqrt p\r).
\end{equation}

Combining \eqref{cot-p} with Lemma \ref{Lem4.3} and \eqref{conjugate}, we obtain
\begin{align*}&\prod_{k=1}^{(p-1)/2}\l(1+\cot\pi\f{ak^2}p\r)
\\=&(1+i)\l((1+i)^2\r)^{(p-3)/4}\f{(-1)^{(p-3)/4}+i}2
\\&\times\f{(-1)^{\f{h(-p)+1}2\cdot\f{p+1}4}}
{(-1)^{\f{h(-p)+1}2}(\f ap)\sqrt p i}\l(s_p+\l(\f ap\r)t_p\sqrt p\r)
\\=&(1+i)(2i)^{(p-3)/4}\f{(-1)^{(p-3)/4}+i}{2i}(-1)^{\f{h(-p)+1}2\cdot\f{p-3}4}\l(\f ap\r)\f{s_p+(\f ap)t_p\sqrt p}{\sqrt p}
\\=&(-1)^{\lfloor\f{p-3}8\rfloor+\f{h(-p)-1}2\cdot\f{p-3}4}2^{(p-3)/4}\l(t_p+\l(\f ap\r)\f{s_p}{\sqrt p}\r).
\end{align*}
This proves \eqref{cot43}.

By Lemma \ref{Lem-4.3},
\begin{align*}\prod_{k=1}^{(p-1)/2}(1+e^{2\pi iak^2/p})=&\prod_{k=1}^{(p-1)/2}\f{1-e^{2\pi i2ak^2/p}}{1-e^{2\pi iak^2/p}}
\\=&\f{(-1)^{\f{h(-p)+1}2}(\f{2a}p)\sqrt p\,i}{(-1)^{\f{h(-p)+1}2}(\f{a}p)\sqrt p\,i}
=\l(\f 2p\r)=(-1)^{(p+1)/4}.
\end{align*}
Combining this with \eqref{tan-p} and \eqref{conjugate}, we get
\begin{align*}
&\prod_{k=1}^{(p-1)/2}\l(1+\tan\pi\f{ak^2}p\r)
\\=&(i-1)((i-1)^2)^{(p-3)/4}(-1)^{\f{h(-p)-1}2\cdot\f{p+1}4}
\\&\times\f{(-1)^{(p-3)/4}+i}2\l(s_p+\l(\f ap\r)t_p\sqrt p\r)
\\=&(-2i)^{(p-3)/4}(-i)^{(3-(\f 2p))/2}(-1)^{\f{h(-p)-1}2\cdot\f{p+1}4}\l(s_p+\l(\f ap\r)t_p\sqrt p\r)
\\=&(-1)^{\lfloor\f{p+1}8\rfloor+\f{h(-p)+1}2\cdot\f{p+1}4}2^{(p-3)/4}\l(s_p+\l(\f ap\r)t_p\sqrt p\r).
\end{align*}
This proves \eqref{tan43}.

In view of the above, we have completed the proof of Theorem \ref{Th1.4}. \qed

\begin{lemma}\label{Lem4.3}
Let $p\eq1\pmod4$ be a prime. Then
\begin{equation}\label{-3p} (-1)^{|\{1\ls k< \f p3:\ (\f kp)=-1\}|}(-3)^{(p-1)/4}\eq\begin{cases}1\pmod p&\t{if}\ p\eq1\pmod{12},
\\\f{p-1}2!\pmod p&\t{if}\ p\eq5\pmod {12}.\end{cases}\end{equation}
\end{lemma}
\Proof. In 1905, Lerch (cf. \cite{HW}) proved that
$$h(-3p)=2\sum_{1\ls k<p/3}\l(\f kp\r),$$
where $h(-3p)$ is the class number of the quadratic field $\Q(\sqrt{-3p})$.
By \cite[Lemma 14]{WC},
$$(-3)^{(p-1)/4}\eq\begin{cases}(-1)^{h(-3p)/4}\pmod p&\t{if}\ p\eq1\pmod{12},
\\(-1)^{(h(-3p)-2)/4}\f{p-1}2!\pmod p&\t{if}\ p\eq5\pmod {12}.\end{cases}$$
Thus, if $p\eq1\pmod{12}$ then
$$(-3)^{(p-1)/4}\eq(-1)^{\f12\sum_{k=1}^{(p-1)/3}((\f kp)-1)+\f{p-1}6}=(-1)^{|\{1\ls k<\f p3:\ (\f kp)=-1\}|}\pmod p;$$
similarly, if $p\eq5\pmod{12}$ then
$$(-3)^{(p-1)/4}\eq(-1)^{|\{1\ls k<\f p3:\ (\f kp)=-1\}|}\,\f{p-1}2!\pmod p.$$
This proves \eqref{-3p}. \qed

\begin{lemma}\label{Lem4.4}
Let $p>3$ be a prime. Then
\begin{equation}\label{relation}
G_p(-\omega)=\begin{cases}(\f{-1}p)\overline{G_p(\omega)}/G_p(\omega)&\t{if}\ p\eq1,3\pmod 8,
\\(\f 3p)\omega^{(\f p3)-1}/(G_p(\omega)\overline{G_p(\omega)})&\t{if}\ p\eq 5,7\pmod8.
\end{cases}
\end{equation}
\end{lemma}
\Proof. Observe that
$$G_p(\omega)G_p(-\omega)=\prod_{k=1}^{(p-1)/2}(e^{2\pi i2k^2/p}-\omega^2)=(-1)^{(p-1)/2}\prod_{k=1}^{(p-1)/2}
\overline{\omega-e^{2\pi i(-2)k^2/p}}.$$
If $(\f{-2}p)=1$, then
$$\prod_{k=1}^{(p-1)/2}(\omega-e^{2\pi i(-2)k^2/p})=G_p(\omega).$$
If $(\f{-2}p)=-1$, then
$$G_p(\omega)\prod_{k=1}^{(p-1)/2}(\omega-e^{2\pi i(-2)k^2/p})=\prod_{r=1}^{p-1}(\omega-e^{2\pi ir/p})=\f{\omega^p-1}{\omega-1}=\l(\f p3\r)\omega^{1-(\f p3)}.$$
Therefore \eqref{relation} holds. \qed

\medskip
\noindent{\it Proof of Theorem 1.5}. As $(\f{-1}p)=1$, we have
\begin{align*}G_p(\omega)^2=&\prod_{k=1}^{(p-1)/2}(\omega-e^{2\pi ik^2/p})(\omega-e^{-2\pi ik^2/p})
\\=&\prod_{k=1}^{(p-1)/2}\l(\omega^2+1-\omega(e^{2\pi i k^2/p}+e^{-2\pi ik^2/p})\r)
\\=&(-\omega)^{(p-1)/2}\prod_{k=1}^{(p-1)/2}(1+e^{2\pi ik^2/p}+e^{-2\pi ik^2/p})
\\=&\omega^{(p-1)/2}\prod_{k=1}^{(p-1)/2}e^{-2\pi ik^2/p}\f{1-e^{2\pi i3k^2/p}}{1-e^{2\pi ik^2/p}}.
\end{align*}
Note that
$$\sum_{k=1}^{(p-1)/2}k^2=\f16\cdot\f{p-1}2\cdot\f{p+1}2\l(2\cdot\f{p-1}2+1\r)=p\f{p^2-1}{24}\eq0\pmod p,$$
and
$$\f{\prod_{k=1}^{(p-1)/2}(1-e^{2\pi i3k^2/p})}{\prod_{k=1}^{(p-1)/2}(1-e^{2\pi i3k^2/p})}
=\f{\sqrt p\ve_p^{-(\f 3p)h(p)}}{\sqrt p\ve_p^{-h(p)}}=\ve_p^{(1-(\f 3p))h(p)}$$
by Lemma \ref{Lem4.2}(i). So we have
$$G_p(\omega)^2=\omega^{(p-1)/2}\ve_p^{(1-(\f 3p))h(p)},$$
and hence
\begin{equation}\label{sign}G_p(\omega)=\da\omega^{(p-1)/4}\ve_p^{(1-(\f p3))h(p)/2}
\end{equation}
for some $\da\in\{\pm1\}$. Write $\ve_p^{h(p)}=a_p+b_p\sqrt p$ with $2a_p,2b_p\in\Z$.
With the help of Lemma 4.2(ii), we have
$$\ve_p^{ph(p)}=(a_s+b_p\sqrt p)^p\eq a_p\eq-\f{p-1}2!\pmod p$$
as in the proof of Theorem 1.3. Thus
\begin{equation}\label{pp}G_p(\omega)^p\eq
\da\omega^{p(p-1)/4}\l(\f p3\r)\l(\f{p-1}2!\r)^{(1-(\f p3))/2}\pmod p.
\end{equation}
On the other hand,
\begin{align*}G_p(\omega)^p=&\prod_{k=1}^{(p-1)/2}(\omega-e^{2\pi ik^2/p})^p
\\\eq&\prod_{k=1}^{(p-1)/2}(\omega^p-1)=\l((\omega^p-1)^2\r)^{(p-1)/4}
\\=&(\omega^{2p}+1-2\omega^p)^{(p-1)/4}=(-3\omega^p)^{(p-1)/4}\pmod p.
\end{align*}
Combining this with \eqref{pp} and Lemma 4.3, we obtain that
\begin{align*}\da\l(\f p3\r)\l(\f{p-1}2!\r)^{(1-(\f p3))/2}\eq&(-3)^{(p-1)/4}
\\\eq&(-1)^{|\{1\ls k<\f p3:\ (\f kp)=-1\}|}\l(\f{p-1}2!\r)^{(1-(\f p3))/2}\pmod p
\end{align*}
and hence
$$\da=(-1)^{|\{1\ls k<\f p3:\ (\f kp)=-1\}|}\l(\f p3\r)=(-1)^{|\{1\ls k\ls\lfloor\f{p+1}3\rfloor:\ (\f kp)=-1\}|}.$$
Combining this with \eqref{sign} we immediately get \eqref{omega}.

In light of \eqref{omega}, we have
\begin{equation}
\label{oomega}G_p(\omega)\overline{G_p(\omega)}=\ve_p^{(1-(\f p3))h(p)}.
\end{equation}
Combining this with \eqref{relation} and \eqref{omega}, we see that if $p\eq1\pmod 8$ then
$$G_p(-\omega)=\f{G_p(\omega)\overline{G_p(\omega)}}{G_p(\omega)^2}
=\f{\ve_p^{(1-(\f p3))h(p)}}{(\omega \ve_p^{h(p)})^{1-(\f p3)}}=\begin{cases}1&\t{if}\ p\eq1\pmod{24},
\\\omega&\t{if}\ p\eq17\pmod{24}.\end{cases}$$
When $p\eq5\pmod 8$, by \eqref{relation} and \eqref{oomega} we have
$$G_p(-\omega)=\l(\f p3\r)(\omega\ve_p^{h(p)})^{(\f p3)-1}=\begin{cases}1&\t{if}\ p\eq13\pmod{24},
\\-\omega\ve_p^{-2h(p)}&\t{if}\ p\eq5\pmod{24}.\end{cases}$$
Therefore \eqref{-omega41} does hold.

In view of the above, we have completed the proof of Theorem \ref{Th1.5}. \qed

\section{Some Conjectures}
\setcounter{lemma}{0}
\setcounter{theorem}{0}
\setcounter{corollary}{0}
\setcounter{remark}{0}
\setcounter{equation}{0}

Let $p>3$ be a prime with $p\equiv3\pmod 8$. Motivated by \eqref{csc} and Theorem 1.3, the author conjectured that
\begin{equation}\label{h}h(-p)=\frac1{2\sqrt p}\sum_{k=1}^{(p-1)/2}\csc2\pi\frac{k^2}p,
\end{equation}
and this conjecture posted to {\tt MathOverflow} (cf. \cite{S19d}) was confirmed by
Prof. Ping Xi.

Inspired by \eqref{sincos0},  for any prime $p\eq7\pmod 8$ and $\da\in\{\pm1\}$ the author guessed that
\begin{equation}\label{4.3}\sum_{k=1}^{(p-1)/2}\f1{1+\da\sin2\pi k^2/p+\cos2\pi k^2/p}=-\f{p+1}4,\end{equation}
and this was confirmed by the author's PhD student Chen Wang who had read the initial version of this paper.

Now we pose some open conjectures on $G_p(x)$ defined in \eqref{S-poly} with $x$ a root of unity.

\begin{conjecture}\label{Conj5.2}  Let $p>3$ be a prime with $p\equiv 3\pmod4$. Then
\begin{equation}\label{S(omega+-)}\begin{aligned}G_p(\omega)=&(-1)^{(h(-p)+1)/2}\left(\frac p3\right)\frac{x_p\sqrt3-y_p\sqrt{p}}2
\\&\times\begin{cases}i&\text{if}\ p\equiv7\pmod{12},
\\(-1)^{|\{1\le k<\frac p3:\ (\frac kp)=1\}|}i\omega&\text{if}\ p\equiv11\pmod{12},
\end{cases}\end{aligned}\end{equation}
and
\begin{equation}\label{S(omega-)}\begin{aligned}G_p(\bar\omega)=&(-1)^{(h(-p)-1)/2}\left(\frac p3\right)\frac{x_p\sqrt3+y_p\sqrt{p}}2
\\&\times\begin{cases}i&\text{if}\ p\equiv7\pmod{12},
\\(-1)^{|\{1\le k<\frac p3:\ (\frac kp)=1\}|}i\bar\omega&\text{if}\ p\equiv11\pmod{12},
\end{cases}\end{aligned}\end{equation}
where $(x_p,y_p)$ is the least positive integer solution to the diophantine equation
$$3x^2+4\left(\frac p3\right)=py^2.$$
\end{conjecture}

{\it Example} 5.1. For the primes $p=79,\,227$,  Conjecture \ref{Conj5.2} predicts that
$$G_{79}(\omega)=i\frac{\sqrt{79}-5\sqrt3}2\ \ \text{and}\ \
G_{227}(\omega)=i\omega(1338106\sqrt3-153829\sqrt{227}).$$

\begin{remark}\label{Rem5.2} Let $p>3$ be a prime with $p\eq3\pmod4$. In light of Lemma \ref{Lem4.4},
Conjecture \ref{Conj5.2} implies that
\begin{align}G_p(-\omega)=\begin{cases}\omega&\t{if}\ p\eq11\pmod{24},
\\1&\t{if}\ p\eq 19\pmod{24},
\\(\f 3p)\omega^{(1+(\f 3p))/2}(x_p\sqrt3+y_p\sqrt p)^2/4&\t{if}\ p\eq 7\pmod{8},
\end{cases}\end{align}
where $x_p$ and $y_p$ are defined as in Conjecture \ref{Conj5.2}.
\end{remark}

 For any prime $p>3$, it is easy to see that
$$G_p(e^{2\pi i/6})=G_p(-\bar{\omega})=\begin{cases}\overline{G_p(-\omega)}&\t{if}\ p\eq1\pmod4,
\\1/\overline{G_p(-\omega)}&\t{if}\ p\eq7\pmod{12},
\\\omega/\overline{G_p(-\omega)}&\t{if}\ p\eq11\pmod{12}.
\end{cases}$$

\begin{conjecture}\label{Conj5.5} Let $p>5$ be a prime and let $\zeta$ be any primitive tenth root
of unity. Then
$$G_p(\zeta)=\begin{cases}(-1)^{|\{1\ls k\ls \f {p+9}{10}:\ (\f kp)=-1\}|}
&\t{if}\ p\eq21\pmod{40},
\\(-1)^{|\{1\ls k\ls\f {p+1}{10}:\ (\f kp)=-1\}|}\zeta^{2}&\t{if}\ p\eq 29\pmod{40}.
\end{cases}$$
\end{conjecture}
\begin{remark} For primes $p>5$ with $p\not\eq21,29\pmod{40}$, we are unable to
find the exact value of $G_p(\zeta)$ with $\zeta$ a primitive tenth root of unity.
\end{remark}

\begin{conjecture}\label{Conj5.3} Let $p>3$ be a prime.

{\rm (i)} If $p\eq13\pmod{24}$, then
\begin{equation}\label{12}G_p(e^{\pm2\pi i/12})=i(-1)^{\f{p-5}8+|\{1\ls k<\f p4:\ (\f kp)=\mp1\}|}(x_p\sqrt3-y_p\sqrt p)
\end{equation}
and
\begin{equation}\label{125}G_p(e^{\pm2\pi i\f 5{12}})=i(-1)^{\f{p-5}8+|\{1\ls k<\f p4:\ (\f kp)=\pm1\}|}(x_p\sqrt3+y_p\sqrt p),
\end{equation}
where $(x_p,y_p)$ is the least positive integer solution to the equation $3x^2+1=py^2$.

{\rm (ii)} When $p\equiv19\pmod{24}$, we may write $p=(4x)^2+3y^2$ with $x,y\in\mathbb Z$, and we have
\begin{equation}\label{1219}G_p(e^{\pm2\pi i/12})=(-1)^{(p-19)/24+x}(1\pm i)\frac{1+\sqrt3}2
\end{equation}
and
\begin{equation}\label{12195}G_p(e^{\pm2\pi i\f5{12}})=(-1)^{(p-19)/24+x}(1\pm i)\frac{1-\sqrt3}2.
\end{equation}

{\rm (iii)} If $p\equiv1\pmod{24}$, then
$$(-1)^{\lfloor \frac{h(-p)}2\rfloor+|\{1\le k<\frac p{12}:\ (\frac kp)=1\}|} e^{\pm2\pi i\f{p-1}{48}}
G_p(e^{\pm2\pi i/12})>0.$$
If $p\equiv7\pmod{24}$, then
$$\pm(-1)^{\lfloor \frac{h(-p)}2\rfloor+|\{1\le k<\frac p{12}:\ (\frac kp)=1\}|} e^{\pm2\pi i\f{p-1}{48}}
G_p(e^{\pm2\pi i/12})>0.$$

{\rm (iv)} If $p\equiv5\pmod{12}$, then
$$(-1)^{\lfloor \frac{h(-p)}2\rfloor+|\{1\le k<\frac p{12}:\ (\frac kp)=1\}|} e^{\pm2\pi i\f{5(p-1)}{48}}G_p(e^{\pm2\pi i/12})>0.$$
If $p\eq11\pmod{24}$, then
$$\pm(-1)^{\lfloor \frac{h(-p)}2\rfloor+|\{1\le k<\frac p{12}:\ (\frac kp)=1\}|} e^{\pm2\pi i\f{5(p-1)}{48}}G_p(e^{\pm2\pi i/12})>0.$$
When $p\equiv 23\pmod{24}$, we have
$$(-1)^{\lfloor \frac{h(-p)}2\rfloor+|\{1\le k<\frac p{12}:\ (\frac kp)=1\}|} e^{\pm2\pi i\f{5(p-1)}{48}}G_p(e^{\pm2\pi i/12})<0.$$
\end{conjecture}
\begin{remark}\label{Rem5.3} The author posted to {\tt MathOverflow} (cf. \cite{S19e}) his conjecture that the equation $3x^2+1=py^2$ has integer solutions for each prime $p\eq13\pmod{24}$, and this was confirmed by the user {\tt GH from MO} via the theory of binary quadratic forms.
For any odd prime $p$, it is known that $h(-p)$ is even or odd according as $p$ is congruent to $1$ or $3$ modulo $4$.
\end{remark}

{\it Example} 5.2. For the prime $p=997\eq13\pmod{24}$, Conjecture \ref{Conj5.3} predicts that
$$G_{997}(e^{2\pi i/12})=-i(318334327\sqrt3-17462102\sqrt{997}).$$

We also have some other conjectures for values of $G_p(x)$ at roots of unity, see Conjectures 13.18-13.20 of \cite[pp.\,275-278]{book}.

\Ack. The author would like to thank the anonymous referee for helpful comments.

\end{document}